\input amssym.def
\input amssym.tex
\magnification 1200
\font\sectionfont=cmbx12 at 12pt
\def\noblackbox{\overfullrule=0pt}
\noblackbox
\def\un{\underline}
\def\bu{\bullet}
\def\pa{\partial}
\def\pr{\prime}
\def\align{\eqalign}
\def\a{\alpha}
\def\b{\beta}
\def\s{\sigma}
\def\o{\omega}
\def\O{\Omega}
\def\inf{\infty}
\def\qed{\vrule height4pt 
width3pt depth2pt}
\def\hook{\hookrightarrow}
\def\BC{{ \Bbb C}}

\def\BZ{{ \Bbb Z}}
\def\la{\lambda}
\def\g{\gamma}

\def\calC{{\cal C}}
\def\calG{{\cal G}}
\def\calU{{\cal U}}

\def\pha{\phantom}
\def\Bi{Birkha\"user}
\def\Ac{Academic Press}
\def\sec#1{\vskip 1.5pc\noindent
{\hbox 
{{\sectionfont #1}}}\vskip 1pc}
\def\longto{\longrightarrow}
\def\mapright#1{\smash
{\mathop{
\longto}\limits^{#1}}}
\def\mapdown#1{\Big\downarrow
\rlap{$\vcenter
{\hbox{$\scriptstyle#1$}}$}}
\def\theo#1#2{\vskip 1pc
  {\bf   Theorem\ #1.} {\it #2}
\vskip 1pc}
\def\conj#1#2{\vskip 1pc
  {\bf   Conjecture\ #1.} {\it #2}
\vskip 1pc}

\def\lem#1#2{\vskip 1pc{\bf 
\ Lemma #1.}
  {\it #2}\vskip 1pc}

\def\prop#1#2{\vskip 1pc{\bf 
 Proposition \ #1.} 
{\it #2}}
\def\defi#1#2{\vskip 1pc 
{\bf Definition \ #1.}
{\it #2}\vskip 1pc }

\def\proof#1{\vskip 1pc {\bf  
Proof.} {\rm #1} 
\vskip 1pc}

\def\NSF9598{\footnote*{This research was supported in part by NSF grants
 DMS-9504522 and DMS-9803593.}}
\centerline{{\bf DIFFERENTIABLE COHOMOLOGY
OF GAUGE GROUPS.}}

\vskip .1 in
\centerline{Jean-Luc Brylinski\NSF9598 }  
\vskip .2 in

\sec{Introduction}

There is a well-known theory of differentiable cohomology
$H^p_{diff}(G,V)$ of a Lie group $G$ with coefficients in a
topological vector space $V$ on which $G$ acts differentiably.
This is developed by Blanc in [{\bf Bl}].
It is very
desirable to have a theory of differentiable cohomology for a (possibly
infinite-dimensional) Lie group
$G$, with coefficients in an arbitrary abelian Lie group $A$,
such that the groups
$H^l_{diff}(G,A)$ have the expected interpretations. For instance,
$H^2_{diff}(G,A)$ should classify the Lie group central extensions of
$G$ by $A$. In this paper we introduce such a theory and study
various differentiable cohomology classes for finite-dimensional
Lie groups and for gauge groups. We are mostly interested in the coefficient
group $A=\BC^*$. In that case, we have the exponential exact sequence
relating the differentiable cohomologies with coefficients $\BZ$, $\BC$ and $\BC^*$.
This allows us easily to compute $H^l_{diff}(G,\BC^*)$ for a compact
Lie group $G$: it is isomorphic to the  cohomology $H^{l+1}(G,\BZ)$.
In the case of gauge groups, we construct
various differentiable cohomology classes, including
the central extension of a loop group as a special case. We also prove
reciprocity laws for  gauge groups of differentiable manifolds
 with boundary embedded in a complex manifold, in the spirit
of the Segal-Witten reciprocity law for loop groups.

The definition of $H^l_{diff}(G,A)$ uses simplicial sheaves. We consider
the classifying space $BG$ as a simplicial manifold, which in simplicial
degree $p$ is equal to $G^p$. Then over each manifold $G^p$
we have the sheaf $\un A$ of smooth $A$-valued functions. These sheaves organize
into a simplicial sheaf $\un A$ over $BG$. We then define $H^l_{diff}(G,A)$
to be the degree $l$ hypercohomology of $BG$ with coefficients in this simplicial
sheaf.

Our motivation is to construct differentiable analogs of the
classes in  the cohomology $H^p(G_{\delta},\BC^*)$ of the
discrete group $G_{\delta}$ constructed by Cheeger and Simons
[{\bf Chee-S}] using geodesic simplices. Similar classes
have been constructed by Beilinson using his Chern classes
in Deligne cohomology. We construct a differentiable cohomology
class in $H^{2p-1}_{diff}(G,\BC^*)$ corresponding to
a characteristic class in $H^{2p}(BG,\BZ)$. In fact,
we construct a more powerful holomorphic class
in the holomorphic group cohomology $H^{2p-1}_{hol}(G_{\BC},\BC^*)$
where $G_{\BC}$ is the complexification of $G$.
We conjecture that these classes map to the classes in [{\bf Chee-S}]
under the natural map from the holomorphic cohomology of $G_{\BC}$
to the cohomology of the discrete group $G_{\delta}$.

In the spirit of secondary characteristic classes,
we construct an extension of these differentiable cohomology
classes involving differential forms of degree $0,1,\cdots,p-1$
on the various stages $G^q$ of the simplicial manifold $BG$.
Again the constructions are done holomorphically on $BG_{\BC}$.
The precise content of the construction is that it yields
a class in the Deligne (hyper)-cohomology of $BG_{\BC}$.

Deligne cohomology is defined using a complex of sheaves,
and the corresponding secondary characteristic classes
are hard to describe explicitly. For the purpose of local computations
at the identity of the group, one would need only to consider
the de Rham part of these classes. The meaning of this
is that Deligne cohomology is approximately defined as a fiber product
of integer-valued cohomology and truncated holomorphic de Rham cohomology
over complex-valued cohomology. Working locally, the integer-valued
cohomology disappears and one is left with truncated de Rham
cohomology. Concretely, given a complex Lie group $G_{\BC}$
such a truncated de Rham class will be represented by a family
$(\o_1,\cdots,\o_p)$, where $\o_j$ is a holomorphic
$(2p-j)$-form over $G^j$.

Our starting point is the theory of Beilinson which says that any characteristic
class in $H^{2p}(G,\BZ)$ for a complex Lie group leads to a so-called
Beilinson characteristic class, which may be viewed as a
holomorphic cohomology class with coefficients, not in $\BC^*$
but in a Deligne complex of sheaves. This can be thought of as
a holomorphic enrichment of a differentiable cohomology class
with $\BC^*$-coefficients. The advantage of using Deligne
cohomology as coefficients is that there are transgression maps
defined for Deligne cohomology, and thus for a closed
oriented manifold $X$ of dimension $k<p$ we can construct
a cohomology class in $H^{2p-1-k}(Map(X,G),\BC^*)$ for the gauge
group $Map(X,G)$. This formalism implies easily that these classes satisfy
a reciprocity law in a holomorphic context (Theorem 3.2). This hopefully clarifies
the meaning of the reciprocity laws in [{\bf Br-ML1}] [{\bf Br-ML2}] [{\bf Br-ML4}].
We note that in [{\bf Br-ML1}] it is incorrectly stated that the
Segal-Witten reciprocity law holds true not only for holomorphic gauge
groups on a Riemann surface with boundary, but even for gauge groups of smooth maps.

We tackle the question of writing down explicit cocycles for all
these cohomology classes. Since this appears at present to be
too difficult a goal for the classes on the whole group,
we localize the question in a neighborhood of the identity.
In this way, the topological part of Deligne cohomology
disappears and we are left with a class in truncated de Rham cohomology.
We conjecture (Conjecture 3.3) that the de Rham cohomology class given by the Beilinson
characteristic class coincides with the class given by the Bott-Shulman-Stasheff
differential forms on $G,G^2,\cdots,G^p$. This conjecture then allows
us to find explicit formulas for local group cocycles
with coefficients in $\BC^*$ (or with enriched coefficients)
in a neighborhood of the identity. The formulas are in the spirit
of [{\bf Br-ML3}].

Interesting phenomena appear when we take the derivative of
a differentiable cocycle to get a Lie algebra cocycle.
In the case of a loop group, we take $p=2$ and $k=1$
and  get exactly the Kac-Moody
$2$-cocycle for a loop algebra $Map(S^1,\frak g)$. For arbitrary $p$, if we take $k=p-1$ we
obtain a well-known Lie algebra $p$-cocycle on $Map(X,\frak g)$ which is due
to Tsygan [{\bf Ts}] and
 Loday-Quillen
[{\bf L-Q}] for $\frak g=\frak g\frak l(n)$ and to Feigin [{\bf Fe}] in general.
This cocycle is the direct higher-dimensional 
generalization of the Kac-Moody cocycle.

I thank Brown University for its hospitality during the Summers of 1997 and
1998, when part of the research and writing were done.

I thank Sasha Beilinson for useful conversations and 
 Daniel Delabre for his many useful comments on a first draft of this paper.

\hfill\eject

\sec{1. Differentiable group cohomology}

For a Lie group $G$ acting smoothly on
a complete topological vector space $M$,
there is a notion of differentiable
cohomology $H^p_{diff}(G,M)$ introduced
and studied by Blanc [{\bf Bl}]. This is equal to the cohomology
of the complex of smooth cochains $C^p(G,M)$.

If instead we consider an abelian Lie group $A$
with a smooth action of $G$, the smooth cochain
complex can still be  defined, but the corresponding
cohomology theory is not fully satisfactory. 
We develop here a formalism for differentiable group
cohomology  which is well-adapted to geometric
applications.

In this paper,  a Lie group $G$ means a paracompact Fr\'echet manifold
$G$ equipped with a group structure such that the product map
and the inverse map are smooth, and there is an everywhere
defined exponential map $exp:\frak g\to G$, where $\frak g$
is the Lie algebra of $G$. Note that this is more restrictive than
the definition in [{\bf P-S}].

The classifying space $BG$
is the simplicial manifold
\hfill\break \vskip .2 in
\hbox{\pha{aaaaaaaa}\hbox{$G\times G\times G $}\pha{aa}\lower.25in
\vbox{ \hbox{$\mapright{d_0}$}
\hbox{$\mapright{d_1}$} 
\hbox{$\mapright{d_2}$}\hbox{$\mapright{d_3}$} }
\pha{aa}\hbox{$G\times G$}
\pha{aa}\lower.1505in
\vbox{\hbox{$\mapright{d_0}$}
\hbox{$\mapright{d_1}$} 
\hbox{$\mapright{d_2}$}}
\hbox{\pha{aa}$G$}
\pha{aa}\lower.075in
\vbox{\hbox{$\mapright{d_0}$}
\hbox{$\mapright{d_1}$}}
\pha{aa}\hbox{pt}
\pha{aaaaaaaaaaaaa}(1-1)}
\par

\vskip 1pc

The face maps are given by the usual formulas.

For a paracompact manifold $X$ and a Lie group $A$, we denote by $\un A_X$
the sheaf of smooth $A$-valued functions on $X$.

Given an abelian Lie group $A$, we can then consider the
simplicial sheaf $\un A$ on $BG$, which consists
in the sheaves $\un A_{G^{ p+1}}$ on each $G^{ p+1}$,
together with the transition morphisms
$$d_j^*\un A_{G^p}\to A_{G^{ p+1}}\eqno(1-2)$$
given by pull-back of smooth $A$-valued functions.

\defi{1.1}{The differentiable cohomology groups
$H^j_{diff}(G,A)$ are the hypercohomology groups
$H^j(BG,\un A)$ of the simplicial sheaf $\un A$ over $BG$.}

This hypercohomology may be computed in two ways.

First of all, the abstract method is to use
 a resolution $\un I^{\bu}$ of the simplicial sheaf $\un A$
over $BG$ by a complex of simplicial sheaves $\un I^q_{\bu}$
such that each sheaf $\un I^q_p$ over $G^p$
is acyclic. Then $H^p_{diff}(G,A)$ is the cohomology
of the double complex $\Gamma(G^{ p},\un I^q_p)$.

Another method is to use \v Cech cohomology. This is somewhat complicated to describe
because we need open coverings of each $G^p$, which are related to each
other via the face and degeneracy maps. More precisely, we need to pick
a family
$\calU^{(p)}=
\{ U^{(p)}_{j\in J_p}\}$
of open coverings of $G^p$, with indexing sets $J_p$, where
the sets $J_p$ form a simplicial set. This means

  (1) for each
$p$ and for each $j\in J_p$ and for $0\leq k\leq p$, there is
a map of sets $d_k:J_p\to J_{p-1}$ such that
$$d_kU^{(p)}_j\subseteq  U^{(p-1)}_{d_k(j)}.$$

(2) for each $p$ and each $0\leq k\leq p$, there is a degeneracy map
$s_k:J_{p-1}\to J_p$ such that 

$$s_kU^{(p-1)}_j\subseteq  U^{(p)}_{s_k(j)}.$$

It is  required  that these maps $d_k$ and $s_k$ satisfy
 the standard relations
among face and degeneracy maps. Such a family of  coverings is
 called good if each
covering $\calU^{(p)}$ of $G^p$ is good (i.e., all non-empty
intersections are contractible). Then for such
 a good covering,
we can form the
 \v Cech double complex
$$\calC^q(\calU^{(p)},\un A)=\oplus_{j_0,\cdots,j_q\in J_p}
~C^{\inf}(U^{(p)}_{j_0\cdots j_q},\un A)\eqno(1-3)$$
where the horizontal differential is the alternating sum of
the pull-back maps 
$$ d_i^*:C^{\inf}(U^{(p)}_{d_i(j_0)\cdots d_i(j_q)},\un A)\to
C^{\inf}(U^{(p-1)}_{j_0\cdots j_q},\un A)$$
The cohomology of this double complex computes the hypercohomology
groups $H^p(BG,\un A)$.

The \v Cech double complex gives rise as usual to  a spectral sequence
$$E_1^{p,q}=H^q(G^p,\un A_{G^p})\Rightarrow H^{p+q}_{diff}
(G,A).\eqno(1-4)$$

We note the following

\lem{1.2}{Given an exact sequence $0\to A\to B\to C\to 0$
of abelian Lie groups, such that the projection
map $B\to C$ has local smooth sections, there is a corresponding long exact
sequence
$$\cdots\to H^{p-1}_{diff}(G,C)\to H^p_{diff}(G,A)\to H^p_{diff}(G,B)
\to H^p_{diff}(G,C)\to\cdots.\eqno(1-5)$$}

\proof{This follows easily from an exact sequence of \v cech double complexes.
\qed}

In case $A$ is a topological vector space, we can describe differentiable group
cohomology with coefficients in $A$ as the cohomology of the
complex of differentiable cochains.

\prop{1.3}{If $A$ is a topological vector space, 
$H^p_{diff}(G,A)$ is the cohomology of the
complex $C^p(G,A)=\{ f:G^p\to A, f~{\rm smooth}~\}$
of smooth $A$-valued cochains.}

\proof{This follows since each $\un A_{G^p}$
is an acyclic sheaf. \qed} 

Therefore in that case our definition of  $H^p_{diff}(G,A)$
coincides with that of Blanc [{\bf Bl}].
In general there is a difference between the differentiable
cohomology $H^p_{diff}(G,A)$ introduced here and the
traditional differentiable cohomology groups $H^p_{class}(G,A)$, which are defined
as the cohomology of the complex $C^p(G,A)$ of smooth $p$-cochains
$G^p\to A$. There always is a map
$$H^p_{class}(G,A)\to H^p_{diff}(G,A)$$
which may be defined as follows:
let $\un I^q_p$ be a resolution of the simplicial sheaf $\un A$ over $BG$,
where each sheaf is acyclic.
There is a map of complexes of simplicial sheaves
$\un A_p\to \un I^{\bu}_p$, which induces a morphism on
the complexes of  global sections. The classical cohomology
$H^p_{class}(G,A)$ is equal to the cohomology of the
complex $\Gamma(G^p,\un A_{G^p})$, and the differentiable
cohomology $H^p_{diff}(G,A)$ is equal to the cohomology of the double complex
$\Gamma(G^p,\un I^q_p)$. This gives the map we announced.

On the other hand, we have:

\lem{1.4}{If $A$ is a discrete abelian group, then
$H^p_{diff}(G,A)$ is equal to the topological cohomology
group $H^p(BG,A)$.}

This leads to the expected for $G$ compact.

\prop{1.5}{Let $G$ be a compact Lie group. Then we have a 
canonical isomorphism
$$H^p_{diff}(G,\BC^*)\simeq H^{p+1}(BG,\BZ(1).$$}

\proof{We use the exponential exact sequence
of coefficient groups:
$$1\to \BZ(1)\to \BC\to{\BC^*}\to 1.$$
together with the vanishing of $H^p_{diff}(G,\BC)$
for $p\geq 1$, proved by Blanc [{\bf Bl}]. \qed}

We now describe the low degree differentiable cohomology
groups.

\lem{1.5}{The group $H^1_{diff}(G,A)$ is the group of smooth
homomorphisms $G\to A$.}

\proof{The spectral sequence
(1-5) gives an exact sequence
$$0\to H^1_{diff}(G,A)\to H^0(G,\un A)
\mapright{d_0^*-d_1^*+d_2^*}
 H^0(G\times G,\un A).$$
The group $H^0(G,\un A)$ is the group of smooth mappings
$\phi:G\to A$; then $\phi$ is a group homomorphism
if and only if $\phi$ is in the kernel of $d_0^*-d_1^*+d_2^*$.
\qed}

\prop{1.6}{The group $H^2_{diff}(G,A)$ is the group of isomorphism
classes of central extensions of Lie groups
$$1\to A\to \tilde G\mapright{\pi} G\to 1$$
such that $\pi$ is a locally trivial smooth principal
$A$-fibration.}

\proof{Given such a central extension, pick a good open covering
$(U_j)_{j\in J}$ of $G$ over which $\pi$ has a smooth section $s_j$.
Then we cover $G\times G$ by the open sets
$U^{(2)}_{j_0j_1j_2}$ defined
as
$$U^{(2)}_{j_0j_1j_2}=d_0^*U_{j_0}\cap
d_1^*U_{j_1}\cap d_2^*U_{j_2}.$$
So $U^{(2)}_{j_0j_1j_2}$ is the set of $(g_0,g_1)$ such
that $g_1\in U_{j_0},g_0g_1\in U_{j_1},g_0\in U_{j_2}$.
This covering is indexed by $J^3$.
Similarly we define open coverings $\calU^{(p)}$
of $G^p$. This allows us to construct a \v Cech double complex.
Because $(U_j)$ is a good covering, this double complex
calculates the differentiable cohomology in degrees $\leq 2$.
We can construct a degree $2$ cocycle in this double complex
as follows: First over $U_{j_0j_1}$ we have
$s_1=s_0g_{01}$, where $g_{01}$ is a smooth function
$U_{j_0j_1}\to A$ (the transition cocycle of the covering).
Next, over $U^{(2)}_{j_0j_1j_2}$ we have the function
$h_{j_0j_1j_2}:U^{(2)}_{j_0j_1j_2}\to A$ defined
by 
$$s_{j_2}(g_0)s_{j_0}(g_1)=s_{j_1}(g_0g_1)h_{j_0j_1j_2}.\eqno(1-6)$$
Then $(g_{j_0j_1},h_{j_0j_1j_2})$ is a $2$-cocycle in the
\v Cech double complex.  It is easy to check that a change in the
choices of the sections $s_i$ will change this $2$-cocycle
 by a coboundary.
Conversely, a degree $2$ cohomology class is
 represented by a $2$-cocycle
$(g_{j_0j_1},h_{j_0j_1j_2})$. The $g_{j_0j_1}$ are the transition
cocycles for a principal $A$-bundle $\tilde G\to G$, equipped with
sections $s_i:U_i\to \tilde G$. There is a unique group structure
on $\tilde G$ compatible with the $A$-action and such that
$$s_{j_2}(g_0)s_{j_0}(g_1)=s_{j_1}(g_0g_1)h_{j_0j_1j_2}.$$
This gives a homomorphism from $H^2_{diff}(G,A)$ to the group
of isomorphism classes of central extensions. 
This is inverse to the map previously constructed. \qed}

Therefore it follows that the image of the map
$H^2_{class}(G,A)\to H^2_{diff}(G,A)$ is comprised of the classes
of central extensions which are trivial as a bundle over $G$.
As is well-known, the universal central extension of a loop group
does not have this property, so it cannot be represented by a class
in $H^2_{class}(G,\BC^*)$.

The description of degree $3$ differentiable cohomology
requires the notion of a {\it gerbe} $\calC$ over a manifold $X$
with band the sheaf $\un A_X$. Then  given a simplicial manifold
$X_{\bu}$ we have the notion of a simplicial gerbe over 
a simplicial manifold $X_{\bu}$,
which was introduced in [{\bf Br-ML1}]. This consists of a gerbe
$\calC$ over $X_1$ with band $\un A_{X_1}$, together with

(1)  an equivalence
$\phi:d_0^*\calC\otimes d_2^*\calC\to d_1^*\calC$ of gerbes with band
$\un A_{X_2}$ over $X_2$;

(2) a natural transformation
$$\psi:d_0^*\phi\otimes d_2^*\phi\to d_1^*\phi\otimes d_3^*\phi\eqno(1-7)$$
between equivalence of gerbes over $X_3$.

The natural transformation $\psi$ must satisfy a cocycle condition.

The structures (1) and (2) become somewhat more concrete in the case of
the simplicial manifold $BG$: (1) can be called a multiplicative
structure on the gerbe $\calC$ over $G$. For instance, on the
 level of the fibers $\calC_g$, which are connected groupoids
in which the automorphism group of any object is identified with $A$, we have
an induced equivalence of groupoids
$$\phi_{g_0g_1}:\calC_{g_0}\otimes \calC_{g_1}\to \calC_{g_0g_1}.
\eqno(1-8)$$
Then one can view $\psi_{g_0g_1g_2}$ as  associativity data for these
equivalences $\phi_{g_0g_1}$.

A simplicial gerbe over $BG$ will be called a multiplicative gerbe over $G$.
Namely, the equivalences $\phi_{g_0g_1}$ are not strictly associative, but only associative
up to the natural transformations $\psi_{g_0g_1g_2}$.

We can then state

\prop{1.7}{The group $H^3_{diff}(G,A)$ identifies with the group of
 equivalence classes of  multiplicative gerbes over $G$ with band $\un A_G$.}

\proof{This is a special case of Theorem 5.7 in [{\bf Br-ML1}], which says that
for  $\un A$ the simplicial sheaf over a simplicial manifold $X_{\bu}$
associated to an abelian Lie group $A$, the
hypercohomology $H^p(X_{\bu},A)$ identifies with the group
of equivalence classes of simplicial gerbes over $X_{\bu}$ with band $\un A$.
\qed}

We next briefly discuss the case of a non-trivial
differentiable $G$-module $A$. This means that there is an action
$\mu:G\times A\to A$ of $G$ on $A$, where the mapping $\mu$ is smooth.
Then we define a simplicial sheaf $\un A$ on $BG$ as follows:
again $A_p$ is the sheaf $\un A_{G^p}$. The face maps
$\nu_i:d_i^*\un A_{G^{p-1}}\to \un A_{G^p}$ are defined as follows:

(1) for $i>0$, $\nu_i$ is the pull-back map $d_i^*$ on
 smooth $A$-valued functions;

(2) for $i=0$, we have
$$\nu_i(f)(g_0,\cdots,g_{p-1})=g_0 f(g_1,\cdots,g_{p-1})\eqno(1-9)$$
using the action of $g_0\in G$ on $A$.
We then  define $H^p_{diff}(G,A)$ to be the hypercohomology
of this simplicial sheaf.

We then have an easy generalization of Lemma 1.5, for an exact sequence
of differentiable $G$-modules such that the map $B\to C$ admits a
local smooth section. There is also the following generalization
of Proposition 1.6:

\prop{1.8}{For any differentiable $G$-module $A$, the group
$H^2(G,A)$ identifies with the group of isomorphism classes
of extensions of Lie groups
$$1\to A\to \tilde G\to G\to 1$$
compatible with the given action of $G$ on $A$, and such that
the mapping $\tilde G\to G$ has local smooth sections.}

We will need crucially the holomorphic version of differentiable cohomology.
For this purpose, $G$ will be a complex-analytic Lie group, $A$ will
be an abelian complex-analytic Lie group, $G$ will act on $A$
 in such a way that the mapping
$G\times A\to A$ is holomorphic. Then we can form the simplicial sheaf
$\un A^{hol}_{\bu}$ over $BG$ such that each $\un A^{hol}_p$ is the
sheaf of germs of holomorphic mappings $G^p\to A$. Then we define
the holomorphic cohomology $H^p_{hol}(G,A)$ to be
 the hypercohomology of the simplicial sheaf
$\un A^{hol}_{\bu}$ over $BG$.

In this paper we will need more general coefficients for
differentiable group cohomology than differentiable $G$-modules.
An important type of coefficients is provided by the smooth
Deligne complex of sheaves $\BZ(k)^{\inf}_D$ over any smooth manifold
$X$:

$$\BZ(k)^{\inf}_D=\BZ(k)\to \un E^0_X\to \cdots\to \un E^{k-1}_X,
\eqno(1-10)$$
where $\un E^l_X$ is the sheaf of germs of smooth complex-valued $l$-forms
over $X$. Then we can organize the complexes of sheaves $\BZ(k)^{\inf}_D$
over the cartesian powers $G^p$ into a simplicial complex of sheaves,
which we will also denote by $\BZ(k)^{\inf}_D$.

\defi{1.9}{The smooth Deligne differentiable cohomology groups
\hfill\break
$H^p_{diff}(X,\BZ(k)^{\inf}_D)$ are the hypercohomology
groups of the simplicial complex of\hfill\break sheaves  $\BZ(k)^{\inf}_D$
over $BG$.}

Finally, for a complex manifold $X$ we have the Deligne complex of 
sheaves $\BZ(k)_D$ over $X$:

$$\BZ(k)_D=\BZ(k)\to \O^0_X\to\cdots\to \O^{k-1}_X,\eqno(1-11)$$
where $\O^k_X$ is the sheaf of holomorphic $k$-forms on $X$.
For a complex Lie group $G$, we can then define the holomorphic Deligne
 cohomology groups
$H^p_{hol}(G,\BZ(k)_D)$ are the hypercohomology
groups of the simplicial complex of sheaves  $\BZ(k)_D$
over $BG$. In other words, $H^p_{hol}(X,\BZ(k)_D)$ is the
Deligne cohomology of the simplicial complex manifold $BG$.

\sec{2. Local cohomology and Lie algebra cohomology.}

The differentiable cohomology $H^p_{diff}(G,A)$
introduced in section 1
involves the sheaf cohomology of the sheaf $\un A_{G^p}$
over $G^p$. In practice classes in $H^p_{diff}(G,A)$ are described
using a \v Cech bicomplex for a family of open coverings of the $G^p$.
If we localize at the origin of the group
the situation becomes simpler. We will therefore introduce
the notion of local differentiable (or simply local)
cohomology $H^p_{loc}(G,A)$. This was developed in the general context of
differentiable groupoids by Weinstein and Xu [{\bf W-X}]. We let
$C^{\inf}_{loc}(G^p,A)$ be the group of germs at $(1,\cdots,1)$ of smooth functions
$G^p\to A$. Then we have a standard complex
$$\cdots\to C^{\inf}_{loc}(G^p,A)\mapright{d_0^*-d_1^*+\cdots+(-1)^pd_p^*}
C^{\inf}_{loc}(G^{p+1},A)\cdots\to.\eqno(2-1)$$

We define the local differentiable cohomology groups
$H^p_{loc}(G,A)$ to be the cohomology groups of this complex.

The main result is

\prop{2.1}{There is a canonical map
$H^p_{diff}(G,A)\to H^p_{loc}(G,A)$.}

We will describe this map concretely. We introduce
a good open covering $\calU^{(p)}$ and consider a \v Cech representative
$$(g^0,\cdots,g^p),$$
where $g^j$ is a \v Cech $(p-j)$-cocycle with coefficients in $\un A_{G^j}$.
 We focus on the last term
for which we observe

\lem{2.2}{If $U$ is an open set in $\calU^{(p)}$
which contains $1$, then the function $g^p_U$ over $U$ satisfies the
cocycle condition
$$\sum_{j=0}^p~(-1)^jd_j^*g^p_U=0$$
in a neighborhood of $1\in G^{p+1}$.}

Now $\calU^{(0)}$ is an open covering
of the set $\{ 1\}$, so has at least  one (non-empty) open set
$U_0=\{ 1\}$. Applying an arbitrary composition of degeneracy
maps, we obtain for each $p$ a distinguished element $i_p$
of $J_p$. The corresponding  open set $V_p=U^{(p)}_{i_p}$ contains
the identity element, and we have
$$d_0(V_p)=d_1(V_p)=\cdots=d_p(V_p).$$
It then follows from Lemma 2.2 that $g^p_U(g_1,\cdots,g_p)$ defines
a local differentiable group cocycle. 
This defines our map from differentiable to local differentiable group
cohomologies.

\underbar{Remarks}

(1) The map $H^p_{diff}(G,A)\to H^p_{loc}(G,A)$ can be described from any open set $U\in
\calU^{(p)}$ containing $1$ which satisfies 
$$d_0(U)=d_1(U)=\cdots=d_p(U).\eqno(2-2).$$ Indeed one can show
that the cohomology class of the corresponding
local cocycle is independent of $U$. 
If $U$ and $V$ satisfy (2-2), then we have:
$$g^p_U-g^p_V=\sum_{j=0}^p~(-1)^jg^{p-1}_{U^{\pr}V^{\pr}}\circ d_j,$$
where $U^{\pr}=d_0(U),V^{\pr}=d_0(V)$.

Thus the difference between $g^p_U$ and $g^p_V$ is a coboundary.

(2) There is a more abstract construction of the localization map,
using the notion of topos.
We associate a very natural
``localized topos''  $X_x$ to a topological space $X$ and a point $x\in X$.
This topos is the $2$-direct limit  of the categories $Sheaves (U)$ 
of sheaves
on open neighborhoods $U$ of $x$. The functor of global sections
is exact, since it coincides with the functor $F\mapsto F_x$, where
$F_x$ is the stalk at $x$. There is a simplicial topos $BG_1$
which in degree $p$ is given by the localized topos $(G^p)_1$.
We can then describe $H^p_{loc}(G,A)$
as the hypercohomology $H^p(BG_1,\un A)$ of the simplicial sheaf
$\un A$ over $BG_1$.

There is a natural map of simplicial topoi $j:BG_1\to BG$, and
the map $H^p_{diff}(G,A)\to H^p_{loc}(G,A)$ is given by the inverse image
$j^*$.

\vskip 1.4pc
There is also a notion of local smooth Deligne cohomology
 which we will have to use.
It will be denoted by $H^p_{loc}(G,\BZ(k)^{\inf}_D)$. It is
 described concretely
as the cohomology of the double complex
$K^{pq}$, where
$$K^{pq}=\left\{\align{&\pha{aa}\BZ(k)~{\rm if}~q=0.\cr
&\lim_{\to\atop
 1\in U\subset G^p}~E^{q-1}(U) ~{\rm if}~1\leq q\leq k\cr
&\pha{aa}0~{\rm  if }~q\geq k+1}\right.\eqno(2-3)$$
A more convenient double complex is the multiplicative version,
where we use the multiplicative version
$$0\to \un\BC^*_X\mapright{d~\log} \un E^1_X
\to\cdots\to \un E^{k-1}_X$$
of the smooth Deligne complex $\BZ(k)^{\inf}_D$ on a manifold
$X$. This leads to the double complex $L^{pq}$, where
$$L^{pq}=\left\{\align{&\lim_{\to\atop 1\in U\subset G^p}~C^{\inf}(U,\BC^*)
~{\rm if}~q=1\cr
&\lim_{\to\atop 1\in U\subset G^p}~E^{q-1}(U) ~{\rm if}~1\leq q\leq k\cr
&\pha{\in}0~{\rm otherwise}}\right.\eqno(2-4)$$

Next there is natural mapping from the local cohomology
$H^p_{loc}(G,A)$ to the Lie algebra cohomology
$H^p(\frak g,A)$, which as usual is defined as the cohomology
of the standard complex
$$\cdots\to C^p(\frak g,A)\to C^{p+1}(\frak g,A)\to\cdots\eqno(2-5)$$
where $C^p(\frak g,A)$ is the space of smooth alternating
multilinear maps $\frak g^p\to A$. The differential $d$
is given by the standard formula.

\prop{2.3}{There is a natural mapping
$\phi:H^p_{loc}(G,A)\to H^p(\frak g,A)$
given on the level of local differentiable cocycles
by the formula
$$\phi(c)(\xi_1,\cdots,\xi_p)=\bigl[{\pa^p\over\pa t_1
\cdots\pa t_p}
\sum_{\sigma\in S_p}~\epsilon(\sigma) c(exp(t_{\sigma(1)}\xi_{\sigma(1)}),
\cdots,exp(t_{\sigma(p)}\xi_{\sigma(p)}))\bigr]_{t_i=0}.$$}

Thus we have a diagram of maps
$$H^p_{diff}(G,A)\to H^p_{loc}(G,A)\to H^p(\frak g,A).$$

We discuss briefly the notions of continuous and measurable
cohomologies. If $A$ is a continuous $G$-module, where
$G$ is a topological group, we define over $BG$ the simplicial
sheaf $\un A_{cont}$ which in degree $p$ is the sheaf of continuous $A$-valued
functions on $G^p$. Then the continuous cohomology groups
$H^p_{cont}(G,A)$ are the hypercohomology groups $H^p(BG,\un A_{cont})$.

If $G$ is a locally compact topological group, then we have the notion
of a measurable subset of $G$, and therefore if $A$ is a Lie group 
we can talk of a measurable function $U\to A$, where $U$ is open in $G$
(or more generally in $G^p$). The measurable cohomology
$H^p_{meas}(G,A)$ is then defined as the hypercohomology of the simplicial
sheaf of measurable $A$-valued functions on the $G^p$. However a simplification
occurs here because of the obvious

\lem{2.5}{The sheaf of germs of measurable $A$-valued functions
on $G^p$ is flasque.}

Therefore measurable cohomology can be computed simply as the
cohomology of the complex of measurable $A$-valued cochains,
which reduces to the standard definition of measurable group
cohomology.

\sec{3. Beilinson Characteristic classes and Bott-Shulman-Stasheff forms.}

Let $X_{\bu}$, $Y_{\bu}$ be simplicial complex manifolds,
and let $G$ be a complex Lie group. We say that a holomorphic
simplicial map $\pi:X{\bu}\to Y_{\bu}$
is a holomorphic principal $G$-bundle is we are given a holomorphic left action
of $G$ on each $X_p$, compatible with the face and degeneracy maps,
so that each $\pi_p:X_p\to Y_p$ is a holomorphic principal
$G$-bundle. The Beilinson theory of characteristic classes in Deligne
cohomology applies to holomorphic principal $G$-bundles
$X_{\bu}\to Y_{\bu}$. We have

\theo{3.1}{\rm [{\bf Be}] [{\bf E}]\it  For any holomorphic $G$-bundle
$X_{\bu}\to Y_{\bu}$ and for any $\kappa\in
H^{2p}(BG,\BZ(p))$, there is a canonical class
$\kappa^{Bei}\in H^{2p}_{hol}(Y_{\bu},\BZ(p)_D)$.}

In case $p=2$, $Y$ is a compact complex manifold, and $G$ is simply-connected,
this class has been described quite concretely in [{\bf Br-ML1}].  There is also
a description in terms of holomorphic gerbes over $G$, which is given in
[{\bf Br-ML1}].

The universal case is that of the classifying space $BG$, which
is a simplicial algebraic manifold. 

Deligne cohomology is closely related to so-called Hodge
cohomology. The simplicial de Rham complex $\O^{\bu}_{Y_{\bu}}$
is filtered by the subcomplexes $F^p\O^{\bu}_{Y_{\bu}}$, which
consist of the complexes of sheaves
$$F^p\O^{\bu}_{Y_n}:0\to\cdots\to 0\to \O^p_{Y_n}\to\cdots
\to \O^{dim(Y_n)}_{Y_n}\eqno(3-1)$$
on $Y_n$. The Hodge cohomology groups are the hypercohomology groups
$H^k(Y_{\bu},F^p\O^{\bu}_{Y_{\bu}})$. We then have an exact sequence
$$\align{\cdots\to H^k(Y_{\bu},\BZ(p)_D)&\to H^k(Y_{\bu},F^p\O^{\bu}_{Y_{\bu}})
\oplus H^k(Y,\BZ(p))\to H^k(Y,\BC)\cr
&\to H^{k+1}(Y_{\bu},\BZ(p)_D)\to\cdots}\eqno(3-2)$$
([{\bf Be}], see also [{\bf E-V}]).

In the case of $BG$, the construction of Bott-Shulman-Stasheff gives
a class in \hfill\break
$H^{2p}(Y_{\bu},F^p\O^{\bu}_{Y_{\bu}})$ (see [{\bf B-S-S}]).
We recall the construction: we start with the invariant polynomial
$P$ on the Lie algebra $\frak g$ corresponding to the
characteristic class $\kappa$. This gives the Chern-Weil
representative $P(\O)$ of the characteristic class $\kappa\in H^{2p}(Y,\BC)$
with respect to a principal $G$-bundle $X\to Y$ equipped with a connection
whose curvature is $\O$. Then one can construct secondary characteristic
classes attached to $m$ connections $D_1,\cdots,D_m$ on the bundle.
For this purpose, one introduces the product $Y\times \Delta_{m-1}$
and the pull-back bundle $X\times \Delta_{m-1}\to Y\times \Delta_{m-1}$.
Let $(t_0,\cdots,t_{m-1})$ denote the barycentric coordinates on $\Delta_{m-1}$,
which satisfy $\sum t_i=1$. Then $D=\sum_{j=0}^m~t_jD_j$ is a connection
for the pull-back bundle. Let $R$ denote the curvature of $D$. One then constructs
the $2p-m+1$-form $\kappa_{sec}(X\to Y,D_1,\cdots,D_m)$ on $Y$ as follows:
$$\kappa_{sec}(X\to Y,D_1,\cdots,D_m)=\int_{\Delta_{m-1}}~P(R)\eqno(3-3)$$
by integrating the $2p$-form $P(R)$ on $Y\times\Delta_{m-1}$
in the $\Delta_{m-1}$-direction.

Now consider the universal $G$-fibration $\pi:EG\to BG$ which in degree
$n$ is given by $\pi_n:G^{n+1}\to G^n$:
$$\pi_n(g_0,\cdots,g_n)=(g_0^{-1}g_1,\cdots,g_i^{-1}g_{i+1},\cdots,g_{n-1}^{-1}g_n)\eqno(3-4)$$

Here $G$ acts on $(EG)_p=G^{p+1}$ by
$$g\cdot(g_0,\cdots,g_p)=(gg_0,\cdots,gg_p).$$
There are $n+1$ sections $\s_0,\cdots,\s_n:G^n\to G^{n+1}$ which are characterized
by the fact that the image of $\s_i$ is the manifold of $n+1$-tuples
$(g_0,\cdots,g_n)$ such that $g_i=1$. Each section $\s_i$ induces a flat connection
$D_i$ on the bundle $G^{n+1}\to G^n$. Then the Bott-Shulman-Stasheff form $\o_n$ on
$G^n$ is the secondary characteristic class $\kappa_{sec}(G^{n+1}\to
G^n,D_0,D_1,\cdots,D_n)$. This is a
$2p-n$-form on  $G^n$. We have the following results:

\theo{3.2}{(Bott-Shulman-Stasheff, see \rm  [{\bf B-S-S}])\it
(1) We have $\o_n=0$ for $n>p$.
\vskip .04 in
(2) The family $(\o_1,\cdots,\o_p)$ is a cycle
in the double complex $\O^m(G^n)$.
\vskip .04 in
(3) For any degeneracy map $s_j:G^{m-1}\hook G^m$, the pull-back
$s_j^*\o_m$ vanishes.}

In particular, $\o=\o_1$ is the bi-invariant closed
form on $G$ representing the transgressive class in $H^{2p-1}(G,\BC)$
corresponding to $\kappa$.

The $p$-form $\o_p$ on $G^p$ is very interesting.
Let $v_i$ be a tangent vector to $G$ and let $(v_i)_j$
be the corresponding tangent vector to $G^p$, which lives in the $j$-th copy
of $G$. if $\tau:\{ 1,\cdots,n\}\to \{ 1,\cdots,n\}$ is a map,
then we can form the expression
$$\o_p((v_1)_{\tau(1)},\cdots,(v_p)_{\tau(p)}).\eqno(3-5)$$
It is easily seen that $\o_p(v_1,\cdots,v_p)=0$
if two of the $(v_i)_{\tau(i)}$'s are tangent vectors to the same
of $G$. This follows readily from property (3) in Theorem 3.2.
Thus the expression (3-5) is
determined by its value in the case where $\tau:\{ 1,\cdots,n\}\to \{ 1,\cdots,n\}$
is a permutation. In fact one sees by direct calculation that we
have:
$$\o_p((v_1)_{\tau(1)},\cdots,(v_p)_{\tau(p)})=\epsilon(\tau)
\o_p((v_1)_1,\cdots,(v_p)_p)\eqno(3-6)$$
and that when we evaluate at the origin in $G^p$,
the resulting multilinear map on $\frak g$ defined by
$$\b(v_1,\cdots,v_p)=\o_1((v_1)_1,\cdots,(v_p)_p)\eqno(3-7)$$
is skew-symmetric.

Because of property (1) in Theorem 3.2 we then can view $(\o_1,\cdots,\o_p)$ as a cocycle in
the simplicial complex of sheaves $n\mapsto F^p\O^{\bu}(G^n)$.
This leads to the natural conjecture:

\conj{3.3}{The class of $(\o_1,\cdots,\o_p)$ in $H^{2p}(BG,F^p\O^{\bu}_{BG})$
is the image of the Beilinson class $\kappa_{Bei}(EG\to BG)$ under the mapping
$H^{2p}(BG,\BZ(p)_D)\to H^{2p}(BG,F^p\O^{\bu}_{BG})$.}

This conjecture will play a key role in the rest of the paper.

\sec{4. Formulas for local group cohomology classes.}

In the previous section we introduced an exact sequence relating Deligne
cohomology  and Hodge cohomology of a simplicial complex manifold.
Equivalently, we have the following exact sequence
$$\cdots\to H^k(Y_{\bu},\BZ(p)_D)\to H^k(Y_{\bu},F^p\O^{\bu}_{Y_{\bu}})
\to H^k(Y,\BC/\BZ(p))\to H^{k+1}(Y_{\bu},\BZ(p)_D)\to\cdots\eqno(4-1)$$
([{\bf Be}], see also [{\bf E-V}]). This shows that the difference
between the two cohomologies is entirely due to the singular cohomology
of $Y_{\bu}$ with $\BC/\BZ(p)$-coefficients. When we go the local Deligne cohomology
$H^p_{loc}(BG,\BZ(k)_D)$ we essentially replace $G$ by the limit
of all the open neighborhoods of $1$, which results in the disappearance
of the higher cohomology groups with $\BC/\BZ(p)$-coefficients. Hence we
obtain:

\prop{4.1}{For $p\geq 1$, 
$H^{2p}_{loc}(BG,\BZ(p)_D)$ is isomorphic to the the degree $2p$ cohomology of
 local truncated de Rham  double complex $F^p\O^{\bu}_{loc}(BG)$, with
$(r,s)$-component equal to
$$\align{&\lim_{\to\atop 1\in U\subseteq G^s}~
\Omega^r(U)~{\rm if}~r\geq p\cr
&\pha{aa}0~{\rm otherwise}.}$$}

We have on any complex manifold $Z$ an exact sequence of complexes
of sheaves:
$$0\to F^p\O^{\bu}_Z\to \O^{\bu}_Z\to [\s_{<p}\O^{\bu}_Z]\to 0,\eqno(4-2)$$
where $\s_{<p}\O^{\bu}_Z$ is the truncated
complex
$\O^0_Z\to \cdots\to \O^{p-1}_Z$.
We have a similar exact sequence of complexes of sheaves
over a simplicial complex manifold $Z_{\bu}$.

This induces a boundary map
$$H^l(Z_{\bu},[\s_{<p}\O^{\bu}_{Z_{\bu}}])\to H^{l+1}(Z_{\bu},F^p\O^{\bu}_{Z_{\bu}}).$$
Let make this more concrete when each $Z_q$ is Stein, so that the above cohomologies
are simply computed by the double complexes of global sections. A class
in $H^l(Z_{\bu},[\s_{<p}\O^{\bu}_{Z_{\bu}}])$ is then  given by 
a family of $l-q$-forms $\o_q$ on $Z_l$, for $ q\geq l-p+1$, satisfying
$$\sum_{j=0}^q(-1)^qd_j^*\o_q=(-1)^qd\o_{q-1}$$

We will write down a formula for  a class $\eta$
in $H^{2p+1}(BG,[\s_{<p}\O^{\bu}_{BG}])$ which we conjecture
to be the image
of the Bott-Shulman-Stasheff class under the boundary map
$$H^{2p}(BG,F^p\O^{\bu}_{loc}(BG))
\to H^{2p+1}(BG,[\s_{<p}\O^{\bu}_{BG}])$$
This class $\eta$ is constructed as follows.
We start by picking a contractible open set $U\subset G$ containing $1$.
Just as in [{\bf Br-ML3}], we construct inductively mappings 
$\s_l:U^l\times\Delta_l\to U$ with the following properties:

(1) $\s_0(pt)=1$;

(2) for all $l$ and for $0\leq j\leq l$, denoting by $\delta_j:
\Delta_{l-1}\hook\Delta_l$ the $j$-th face map, we have

$$\s_l(g_1,\cdots,g_l;\delta_j(t_0,\cdots,t_{l-1}))=\biggl\{
\align{&\s_{l-1}(d_j(g_1,\cdots,g_l);t_0,\cdots,t_{l-1})~{\rm if}~j\geq 1\cr
&g_1\cdot \s_{l-1}(g_2,\cdots,g_l;t_0,\cdots,t_{l-1})~{\rm if}~j=0}.
\eqno(4-3)$$

For fixed $(g_1,\cdots,g_l)\in U^l$, the resulting map
$\Delta_l\to U$ will be denoted by $\s_{g_1,\cdots,g_l}$. It is a singular
simplex in $U$ with vertices $(1,g_1,\cdots,g_l)$.

Then  we have  mappings
$$f_{m,q}:U^{m+q-1}\times\Delta_q\to G^m\eqno(4-4)$$
given by  
$$f_{m,q}(g_1,\cdots,g_{m+q-1};t_0,\cdots,t_q)
=(g_1,\cdots,g_{m-1},\s_q(g_m,\cdots,g_{m+q-1};t_0,\cdots,t_q))\eqno(4-5)$$
 Now, $\o_m$ is a $(2p-m)$-form on $G^m$; we can pull it back under $f_{m,q}$
and then integrate it over $\Delta_q$ to get a $(2p-m-q)$-form
$\int_{p_1}~f_{m,q}^*\o_m$ over $U^{m+q-1}$. This $(2p-m-q)$-form will be
denoted by $\b_{m,q}$.

We are now ready to introduce an $l$-form $\eta_l$ over $U^{2p+1-l}$,
defined by
$$\eta_l=\sum_{m+q=2p-l\atop m\geq 1}~\b_{m,q}.\eqno(4-6)$$

Note $b_{m,q}$ vanishes unless $m\leq p$.

We note that $\eta_0=\beta_{1,2p-1}$. Indeed, for $m+q=2p$ and
$m\geq 2$, the function
$\beta_{m,q}$ vanishes, as  the first component of $f_{m,q}$ is 
equal to $g_1$ and thus, when we fix $g_1$ to be a constant, the restriction
of the differential form $\int_{p_1}~f_{m,q}^*\o_m$ to
$\{ g_1\}\times U^{m+q-1}\times \Delta^q$ vanishes.
Now the function $\beta_{1,2p-1}$ is exactly the function
which Cheeger and Simons [{\bf C-S}] use to define a $(2p-1)$ group cocycle of the
discrete group
$G^{\delta}$ with values in $\BC^*$.

Then we can present our  conjecture

\conj{4.2}{The family $\eta_l$ of $l$-forms over $U^{2p+1-l}$ is a cocycle
in the truncated de Rham double complex $\O^r(U^s)$, which represents the image of the
Bott-Shulman-Stasheff class under
the boundary map
$$H^{2p}(BG,F^p\O^{\bu}_{loc}(BG))
\to H^{2p+1}(BG,[\s_{<p}\O^{\bu}_{BG}])\eqno(4-7)$$}

We illustrate this  for $p=2$. The $1$-form $\eta_1$ over $U^2$
 is equal to $\eta_1=\b_{1,2}+\b_{2,1}$. These $1$-forms can be written
down as follows. We introduce the following notations.
Let $\xi_1,\cdots,\xi_{m+q-1}$ denote vector fields on $G$. View
each $\xi_j$ as a tangent vector on the $j$-th factor of $U^{m+q-1}$.
Then $df_{m,q}(\xi_j)$ can be viewed as a section of the pull-back
$f_{m,q}^*TG^m$ of the tangent bundle of $G^m$ to $U^{m+q-1}\times
\Delta_q$. Restricting to a point $(g_1,\cdots,g_{m+q-1})$ of $U^{m+q-1}$,
we may view $df_{m,q}(\xi_j)$ as a section of 
$\s_{g_1,\cdots,g_{m+q-1}}^*TG^m$. 

Our  formula for $\b_{m,q}$ (for $m+q=3$) is then
$$(\b_{m,q})_{(g_1,g_2)}(\xi)
=\int_{\Delta_q}~df_{m,q}(\xi)\rfloor
\s_{g_1,g_2}^*\o_m.\eqno(4-8)$$

For example, $\b_{1,2}$ involves integrating the $2$-form
$df_{1,2}(\xi)\rfloor\o_1$ over the $2$-simplex,
and $\b_{2,1}$ involves integrating
the $1$-form 
$df_{2,1}(\xi)\rfloor\o_2$
over the $1$-simplex.

It should be noted that the presence of $\b_{2,1}$ is necessary
to obtain a local holomorphic  cocycle.

We now prove Conjecture 4.2 in the case $p=2$.
In the proof we will use 3 types of facts:

(1) the (relative) Stokes theorem in the case of a projection
$p_1:X\times\Delta^q\to X$. For $\b$ a differential form
on $X\times\Delta^q$, this gives

$$d\int_{p_1}\b=\int_{p_1}d\b+\sum_{j=0}^q\int_{p'_1}(Id\times\delta_j)^*\b,$$
where $p'_1: X\times\Delta^{q-1}\to X$ is the projection.

(2) the fiber square principle which says that for a fiber square
$\matrix{T&\mapright{g}&X\cr \mapdown{q}&&\mapdown{p}\cr
Z&\mapright{f}&Y}$ in which $p$ is a smooth proper fibration,
and for a differential form
$\b$ on $X$, we have 
$$f^*\int_p\b=\int_qg^*\b.$$

(3) the cocycle relation
$$\sum_{j=0}^3(-1)^jd_j^*\o_2=0.$$

(4) the cocycle relation
$$\sum_{j=0}^2(-1)^jd_j^*\o=d\o_2.$$

First we replace the class $(\o,\o_2)$ in the hypercohomology
of $F^2\O^{\bu}_{loc}(BG)$ by a cohomologous cocycle.
More precisely we add to it the total coboundary of $\beta_{1,1}\in \O^2_{loc}(G)$.
This coboundary has components $-d\b_{1,1}$ and $\sum_{j=0}^2(-1)^j
d_j^*\b_{1,1}$. We compute $d_j^*\b_{1,1}$ using the 
fiber square principle. We find
$$d_j^*\b_{1,1}=\int_{p_1}\nu_j^*\o,$$
where $\nu_j=f_{1,1}\circ (d_j\times Id):U^2\times \Delta^1\to G$ so that
$$\nu_0(g,1,g_2,t)=\gamma_{g_2}(t)$$

$$\nu_1(g,1,g_2,t)=\gamma_{g_1g_2}(t)$$

$$\nu_2(g,1,g_2,t)=\gamma_{g_1}(t)$$

We want to compare this alternating sum with $d\eta_1$. We first compute
$d\b_{2,1}$. Using the Stokes theorem and the cocycle relation
we have
$$\align{d\b_{2,1}&=\int_{p_1}f_{2,1}^*d\o_2+\o_2
=\sum_{j=0}^2(-1)^j\int_{p_1}f_{2,1}^*d_j^*\o+\o_2\cr
&=\sum_{j=0}^2(-1)^j\int_{p_1}\a_j^*\o+\o_2}$$
where $\a_j=d_jf_{2,1}:U^2\times \Delta^1\to U$ so that

$$\a_0(g_1,g_2,t)=\gamma_{g_2}(t)$$

 $$\a_1(g_1,g_2,t)=g_1\cdot \gamma_{g_2}(t)$$

$$\a_2(g_1,g_2,t)=g_1$$

Now using Stokes' theorem we find

$$d\b_{1,2}=\sum_{j=0}^2(-1)^j\int_{p_1}\gamma_j^*\o,$$
where $\a_j=f_{1,2}\circ (Id\times \delta_j):U^2\times\Delta^1\to U$
is given by 

$$\gamma_0(g_1,g_2,t)=g_1\cdot \gamma_{g_2}(t)$$

$$\gamma_1(g_1,g_2,t)=\gamma_{g_1g_2}(t)$$

$$\gamma_2(g_1,g_2,t)=\gamma_{g_1}(t)$$

We then obtain

$$\sum_{j=0}^2(-1)^jd_j^*\b_{1,1}=d\b_{1,2}+d\b_{2,1}-\o_2=d\eta_1-\o_2.$$

We conclude that the Bott-Shulman-Stasheff cocycle $(\o,\o_2)$
is cohomologous to $(0,d\eta_1)$. In order to prove the conjecture for $p=2$,
all we have left to do is to show that $(\eta_0,\eta_1)$ is indeed a cocycle,
i.e.,

$$dh=d_0^*\eta_1-d_1^*\eta_1+d_2^*\eta_1-d_3^*\eta_1\eqno(4-10)$$

We first compute
 $\sum_{j=0}^3(-1)^jd_j^*\b_{2,1}$. By the fiber square principle
we have
$d_j^*\b_{2,1}=\int_{p_1}\phi_j^*\o_2$, where $p_1:U^2\times \Delta_1\to U^2$
is the projection and the map $\phi_j=f_{2,1}\circ (d_j\times Id):U^3\times \Delta_1
\to U^2$ is given by

$$\phi_0(g_1,g_2,g_3,t)=(g_2,\gamma_{g_3}(t));$$

$$\phi_1(g_1,g_2,g_3,t)=(g_1g_2,\gamma_{g_3}(t));$$

$$\phi_2(g_1,g_2,g_3,t)=(g_1,\gamma_{g_2g_3}(t));$$

$$\phi_3(g_1,g_2,g_3,t)=(g_1,\gamma_{g_2}(t)).$$

Here $\gamma_g(t)=\sigma_1(g,t)$.

Similarly we have $d_j^*\b_{1,2}=\int_{p_1}\psi_j^*\o$, where
$\psi_j:U^3\times\Delta_2\to U$ is given by

$$\psi_0(g_1,g_2,g_3,y)=\s_{g_2,g_3}(y)$$

$$\psi_1(g_1,g_2,g_3,y)=\s_{g_1g_2,g_3}(y)$$

$$\psi_2(g_1,g_2,g_3,y)=\s_{g_1,g_2g_3}(y)$$

$$\psi_2(g_1,g_2,g_3,y)=\s_{g_1,g_2}(y),$$
where $y$ denotes a point in $\Delta_2$.

Next we find from Stokes'  theorem that

$$dh=A-d_1^*\b_{1,2}+d_2^*\b_{1,2}-d_3^*\b_{1,2}$$
where $A=\int_{p_1}\la^*\o$, for the mapping
$\la:U^3\times\Delta_2\to U$ defined by
$$\la(g_1,g_2,g_3,y)=g_1\cdot \s_{g_2,g_3}(y).$$
So we are led to consider the mapping
$f_{2,2}:U^3\times\Delta_2\to U^2$ which satisfies
$$d_0\circ f_{2,2}=\psi_0,d_1\circ f_{2,2}=\la,d_2\circ f_{2,2}=p_1.$$

Then we get from the cocycle relation (4):

$$\int_{p_1}df_{2,2}^*\o_2=\sum_{j=0}^2(-1)^j\int_{p_1}f_{2,2}^*d_j^*\o
=d_0^*\b_{1,2}-A+\int_{p_1}p_1^*\o.$$
This is evaluated by Stokes' theorem to be equal
to $d_2^*\b_{2,1}-d_3^*\b_{2,1}+\int_{p_1}\theta^*\o_2$, where
$\theta(g_1,g_2,g_3,t)=(g_1,g_2\cdot \g_{g_3}(t))$.
We now claim that
$$d_0^*\b_{2,1}-d_1^*\b_{2,1}+\int_{p_1}\theta^*\o_2=0.$$

To prove this, we note that $\sum_{j=0}^3(-1)^j\int_{p_1}f_{2,2}^*d_j^*\o_2=0$
by the cocycle relation. It is easy to see that the term corresponding
to $j=3$ vanishes, and that the other three terms are
$d_0^*\b_{2,1}$, $d_1^*\b_{2,1}$ and $\int_{p_1}\theta^*\o_2$.

We finally gather all this information to obtain:

$$\align{dh&=A-d_1^*\b_{1,2}-d_2^*\b_{1,2}+d_3^*\b_{1,2}\cr
&=\sum_{j=0}^3(-1)^jd_j^*\b_{1,2}+d_2^*\b_{2,1}-d_3^*\b_{2,1}
+\int_{p_1}\theta^*\o_2\cr
&=\sum_{j=0}^3(-1)^jd_j^*\b_{1,2}+\sum_{j=0}^3(-1)^jd_j^*\b_{2,1}}$$
which proves the statement.

Conjecture 4.2 implies  that the family
$(\eta_l)$ will represent the class
of $\kappa$ in the local (differentiable) Deligne cohomology $H^{2p}_{loc}
(G,\BZ(p)_D)$,
assuming the validity of Conjecture 3.3.
As evidence for Conjecture 3.3, we note that the 
 formula it implies for the local Deligne cohomology class in the case $p=2$
is very similar to the formula
in [{\bf Br-ML1}]. 
\hfill\eject

\sec{5. Differentiable cohomology classes for gauge groups.}

In this section, we investigate the differentiable cohomology classes
for gauge groups $Map(X,G)$ which result from transgression
of a characteristic class
$\kappa\in H^{2p}(BG,\BZ(2p))$. As a motivation, we briefly
recall from [{\bf Br-ML1}] the case of central extensions
of a smooth loop group $LG=Map(S^1,G)$.
 In this case one starts from a characteristic class
$\kappa\in H^{4}(BG,\BZ(2))$, and by transgression
in the diagram of simplicial manifolds
$$\matrix{(LG)^n\times S^1&\mapright{ev}&G^n\cr
\mapdown{p_1}&&\cr
(LG)^n&&}\eqno(5-1)$$
one obtains a class in $H^3(B(LG),\BZ(2)^{\inf}_D)$, hence in particular
a class in \hfill\break
$H^3(B(LG),\BZ(1)^{\inf}_D)$. Since the complex of sheaves $\BZ(1)^{\inf}_D
$ on a manifold $X$ is quasi-isomorphic to $\un\BC^*_X[1]$, this can be viewed
as a class in $H^2(L(BG),\un\BC^*)=H^2_{diff}(LG,\BC^*)$.
According to Proposition 1.6, this gives a central extension of Lie groups
$$1\to \BC^*\to \widetilde{LG}\to LG\to 1.\eqno(5-2)$$
A different approach to the central extensions is given in [{\bf Pi}].

In fact it is more rewarding to consider the holomorphic
 analog of this construction, for a complex Lie group $G$.
The construction is then performed using Beilinson
characteristic classes. Then the result is a holomorphic central extension.
This is the point of view developed in [{\bf Br-ML1}].

In this case, assuming Conjecture 3.4, we can write down at least a local differentiable
group cocycle for the central extension, using the $1$-form
$\eta_1=\b_{1,2}+\b_{2,1}$ over $U^2$. We transgress this $1$-form in the
evaluation diagram
$$\matrix{LU^2\times I&\mapright{ev}&U^2\cr
\mapdown{p_1}&&\cr
LU^2&&}$$
to obtain a smooth function
$c=\int_{p_1}ev^*\eta_1$ on $U^2$.

This function is given by the following formula
$$c(g_1,g_2)=a(g_1,g_2)+b(g_1,g_2),\eqno(5-3)$$
where

(1) $a(g_1,g_2)=\int_{S^1\times \Delta_2}~\phi^*\o$,
for the mapping $\phi:S^1\times \Delta_2\to G$ given by
$$\phi(u;t_0,t_1,t_2)=f_{1,2}(g_1(u),g_2(u);t_0,t_1,t_2)
=\sigma_{g_1(u),g_2(u)}(t_0,t_1,t_2)
;\eqno(5-4)$$

(2) $b(g_1,g_2)=\int_{S^1\times \Delta_1}~\psi^*\o_2$,
for the mapping $\psi:S^1\times\Delta_1\to G^2$ described by
$$\psi(u;t_0,t_1)=f_{2,1}(g_1(u),g_2(u):t_0,t_1)
=(g_1(u),\s_{g_2(u)}(t_0,t_1)).\eqno(5-5)$$

It is interesting to compute the corresponding Lie algebra
$2$-cocycle, obtained as in \S 2 by differentiating
$c$ at $(1,1)\in U^2$ and skew-symmetrizing. Let $\xi_1,\xi_2$
belong to the Lie algebra $L\frak g$ of $LG$. We observe
that the coefficient of $y_1y_2$ in $a(exp(y_1\xi_1),exp(y_2\xi_2))$
is $0$, because the $3$-dimensional integral is clearly
$O((|y_1|+|y_2|)^3)$.

To evaluate the dominant term of $b(exp(y_1\xi_1),exp(y_2\xi_2))$, we may
 as well assume that
$$f_{2,1}(exp(y_1\xi_1(u)),exp(y_2\xi_2(u));t_0,t_1)
=(exp(y_1\xi_1(u)),exp(t_0y_2\xi_2(u)),$$
so that
$$\psi(u;t_0,1-t_0)=(exp(y_1\xi_1(u)),exp(t_0y_2\xi_2(u)).$$
The derivative of this mapping is given (up to terms which are
$O((|y_1|+|y_2|)^3)$), by
$${\pa\over \pa u}\mapsto (y_1{d\xi_1(u)\over du},
t_0y_2{d\xi_2(u)\over du})\pha{a},
\pha{a}
{\pa\over \pa t_0}\mapsto (0,y_2\xi_2(u)).$$
Now using the expression
$$\o_2=3k~Tr(g_1^{-1}dg_1~dg_2g_2^{-1})$$
we find
$${\pa^2\over \pa y_1\pa y_2}
b(exp(y_1\xi_1),exp(y_2\xi_2))=3k~
\int_0^1~Tr({d\xi_1\over du}\xi_2(u)) du.$$
As this expression is clearly already skew-symmetric in $\xi_1,\xi_2$
its skew-symmetrization gives the $2$-cocycle
$$\a(\xi_1,\xi_2)=6k~
\int_0^1~Tr({d\xi_1\over du}\xi_2(u)) du.\eqno(5-6)$$
This is the standard Kac-Moody cocycle [{\bf Ka}] [{\bf Mo}].

The generalization of these constructions
that we will discuss involves any closed oriented
$d$-dimensional manifold $X$. Then we have the smooth
 gauge group
$Map(X,G)$ comprised of the smooth maps $X\to G$. For $X$ compact
(possibly with boundary), this is
a Lie group with Lie algebra $Map(X,\frak g)$, which 
is a Fr\'echet vector space.

\theo{5.1}{For any characteristic class
$\kappa\in H^{2p}(BG,\BZ(p))$ and any
closed oriented
 manifold $X$ of dimension $k\leq p$, there results a class
in differentiable Deligne cohomology
 $\kappa(X)\in H^{2p-k}(B~Map(X,G),\BZ(p-k)^{\inf}_D)$.}

This class is simply obtained using transgression
in smooth Deligne cohomology for the evaluation
diagram
$$\matrix{Map(X,G)^p\times X&\mapright{ev}&X^p\cr
\mapdown{p_1}&&\cr
Map(X,G)^p&&}\eqno(5-7)$$

We then have a geometric reciprocity law.
This is phrased in terms of a real \hfill\break $(k+1)$-submanifold
with boundary inside a complex manifold of dimension $k$.
The case $k=1$ gives the Segal-Witten reciprocity law
for Riemann surfaces with boundary.
We note that in [{\bf Br-ML1}] it is incorrectly stated that the
reciprocity laws for groups of smooth group-valued maps as opposed
to groups of holomorphic maps.

\theo{5.2}{Let $M$ be a smooth compact manifold of complex dimension $k$
and let $X\subset M$ be a real $(k+1)$-dimensional submanifold
with boundary $\pa X$. Let $G$ be a complex Lie group and let $\calG$ be the Lie group
consisting of smooth maps
$X\to G$  which have a holomorphic extension to some neighborhood of $X$.
Then the pull-back of
$\kappa(\pa X)$ to
$H^{2p-k}(B~\calG,\BZ(p-k)_D^{\inf})$
is trivial.}

In the case of Riemann surfaces with boundary,  the reciprocity law is
more precise than the above theorem: it says not only that
 the pull-back of some central extension splits, but even
that it has a canonical splitting. 

Then we have

\theo{5.3}{Under the assumptions of Theorem 5.2, the class
$\kappa(\pa X)$\hfil\break
$\in H^{2p-k}(B~\calG,\BZ(p-k)_D)$
lifts canonically to a class
in the relative differentiable Deligne cohomology group
$$H^{2p-k}(B~\calG\to B~Map(\pa X,G),\BZ(p-k)_D^{\inf}).$$}

This is the analog in differentiable group  cohomology of a theorem proved
for group cohomology in [{\bf Br-ML2}].

We will write down an explicit formula for the class
in local Deligne cohomology $H^{2p-k}_{loc}(Map(X,G),\BZ(p-k)^{\inf}_D)$
associated to the characteristic class $\kappa$.
Recall from section 4 the cocycle $(\eta_1,\eta_2,\cdots,\eta_p)$
in the double complex $\Omega^q(G^p)$ associated to $\kappa$.
Here $\eta_l$ is an $l$-form over $G^{2p-1-l}$.
Then we have
\prop{5.4}{ Assume the validity of Conjecture 3.3.
For a smooth closed  manifold $X$ of dimension $k$, and for a characteristic class
$\kappa\in H^{2p}(G,\BZ(p))$, the corresponding class
in local Deligne cohomology $H^{2p-k}_{loc}(Map(X,G),\BZ(p-k)^{\inf}_D)$
is represented by the family $\int_X\eta_l$ of $(l-k)$-forms
over $Map(X,G)^{2p-1-l}$, where $\int_X\eta_l$ is the transgression
of $\eta_l$ in the evaluation diagram
$$\matrix{Map(X,G)^p\times X&\mapright{ev}&X^p\cr
\mapdown{p_1}&&\cr
Map(X,G)^p&&}\eqno(5-8)$$}

We now  obtain a Lie
algebra cocycle simply by differentiating at the origin of $G$.
Introduce a mapping
$$\phi_{m,q}:Map(X,G)^{2p-k-1}\times X\times\Delta_q\to G^{2p-k-1}$$
by
$$\phi_{m,q}(g_1,\cdots,g_{2p-k-1};x;(t_0,\cdots,t_q))=
f_{m,q}(g_1(x),\cdots,g_{2p-k-1}(x);(t_0,\cdots,t_q))$$
For given $(g_1,\cdots,g_{2p-k-1})\in Map(X,G)^{2p-k-1}$, denote
by $[\phi_{m,q}]_(g_1,\cdots,g_{2p-k-1})$ the corresponding
mapping $ X\times\Delta_q\to G^{2p-k-1}$.

Then we see  easily that
the Lie algebra cocycle defined from differentiating
the class of Theorem 5.1 
is obtained by skew-symmetrizing the differential
at the origin of $Map(X,G)^{2p-k-1}$ of the local functional
$$(g_1,\cdots,g_{2p-k-1})\mapsto \sum_{m+q=2p-k}
c_{m,q}(g_1,\cdots,g_{2p-m-1}),$$
where $m$ ranges over $1,\cdots,p$.
Here we have put
$$c_{m,q}(g_1,\cdots,g_{2p-k-1})=\int_{X\Delta_q}~
[\phi_{m,q}]_{(g_1,\cdots,g_{2p-m-1})}^*\o_m.$$
We can find an explicit formula for the
partial derivative

$${\pa^{2p-k-1} c_{m,q}(exp(t_1\xi_1),\cdots,exp(t_{2p-m-1}\xi_{2p-k-1}))
\over \pa t_1\cdots \pa t_{2p-k-1}}$$
evaluated at $t_j=0$
 as follows.
First of all, we see easily that this partial derivative will
be $0$ unless $m=k$. In that case,
introduce the value $(\o_m)_1$ of $\o_m$ at the identity, which
is a $2p-m-1=2p-k-1$-multilinear form on $\frak g$. Given
$(\xi_1,\cdots,\xi_{2p-k-1})$ in $Map(X,\frak g)$, we can write down
a $k$-form $\a$ on $X$, whose expression in terms of local coordinates
$(x_1,\cdots,x_k)$ is
$$\a=(k!)^{-1}
[\o_m]_1(\xi_1,\cdots,\xi_{2p-2k-1},{\pa\xi_{2p-2k}\over \pa x_1},\cdots,
{\pa\xi_{2p-m-1}\over \pa x_k})~dx_1\wedge\cdots\wedge dx_k.\eqno(5-9)$$

Then we have

\prop{5.5}{The Lie algebra $(2p-k-1)$-cocycle obtained by differentiating the class in
$H^{2p-k-1}(Map(X,G),\BC^*)$ is equal to the skew-symmetrization of the cochain 
$(\xi_1,\cdots,\xi_{2p-2k-1})\mapsto\int_X\a$ where $\a$ is as in (5-9).}

Consider for instance the case $k=p-1$. In  that case we have
$m=k=p-1$ hence  $q=2p-k-p+1=1$. We then get the following
Lie algebra cocycle:

$$c(\xi_1,\cdots,\xi_p)=\int_X \o_p(\xi_1,d\xi_2,\cdots,d\xi_p)\eqno(5-10).$$
This is a direct generalization of the Kac-Moody cocycle which was given by Feigin
[{\bf Fe}]. We refer the reader to Teleman [{\bf Te}] for a detailed discussion
of this and more general Lie algebra homology classes.

For $\frak g=\frak g\frak l(n$, the above 
construction is consonant with
the results of Tsygan [{\bf Ts}] and of Loday-Quillen [{\bf L-Q}]
on the Lie algebra of infinite matrices over a $\BC$-algebra
$A$. Recall that the Lie algebra homology
$$H_{\bu}(M_{\inf}(A),\BC)=\lim_{\to\atop n}~H_{\bu}(M_n(A),\BC)$$
is a Hopf algebra, and that Tsygan and Loday-Quillen prove that
the primitive part $Prim~H_{\bu}(M_{\inf}(A),\BC)$
of this Hopf algebra identifies with the cyclic homology of $A$.
More precisely, there is a shift of $1$ in the degrees so that
$$Prim_pH_{\bu}(M_{\inf}(A),\BC)\tilde{\to}
HC_{p-1}(A).$$
This applies to $A=C^{\inf}(X)$, viewed as a Fr\'echet algebra;
 the cyclic homology of $C^{\inf}(X)$ is defined taking this topology onto account.
Then the cyclic homology of $C^{\inf}(X)$, computed by Connes [{\bf Co}],
Tsygan [{\bf Ts}], Loday-Quillen [{\bf L-Q}] is given by
$$HC_i(C^{\inf}(X))=\biggl\{\align{&[E^i(X)/d E^{i-1}(X)]\oplus H^{i-2}(X)
\oplus\cdots~{\rm if}~i<k\cr
 &\oplus_{n\in\BZ}~H^{i+2n}(X)~{\rm if}~i\geq k.}$$
Dually, the primitive part in degree $p$ of the Lie algebra cohomology
of $M_{\inf}(A)$ identifies with the cyclic cohomology group
$HC^{p-1}(A)$, which is computed by Connes:
$$HC^i(C^{\inf}(X))=\bigl\{\align{&T_i(X)_{cl}\oplus 
H_{i-2}(X)
\oplus\cdots~{\rm if}~i<k\cr
 &\oplus_{n\in\BZ}~H_{i+2n}(X)~{\rm if}~i\geq k.}$$
Here $T_i(X)_{cl}$ denotes the space of closed 
degree $i$ currents on $X$.
We can also consider the relative cyclic cohomology
$HC^i(C^{\inf}(X),\BC)$ which is obtained by deleting the factor
$H_0(X)$ from $HC_i(X)$ which is present for $i$ even.

So the fundamental class in $H_k(X)$ yields a class
in $HC^k(C^{\inf}(X),\BC)$ which corresponds to a primitive class
in $H^{k+1}(M_{\inf}(C^{\inf}(X))$. This class is precisely
that defined by the Lie algebra cocycle (5-10).

\vskip .2 in

REFERENCES
\vskip .14 in

[{\bf Be1}] A. A. Beilinson, \it
Higher regulators and values of 
L-functions, \rm J. Soviet Math. 
\bf 30 \rm (1985), 2036-2070

[{\bf Bl}] P. Blanc, {\it Cohomologie diff\'erentiable et
changement de groupes}, Ast\'erisque vol. 124-125 (1985), pp. 113-130.

[{\bf B-S-S}] R. Bott, H. Shulman and J. Stasheff, \it
On the de Rham theory of certain classifying spaces,
\rm Adv. in Math. \bf 20 \rm (1976), 535-555

[{\bf Br-ML1}] J-L Brylinski and D. D. McLaughlin, \it 
 The geometry of the first
 Pontryagin class and of line bundles on loop
spaces I, \rm 
 Duke Math. Jour. \bf 35 \rm (1994), 603-638

[{\bf Br-ML2}] J-L Brylinski and D. D. McLaughlin, \it 
 Characteristic classes and
multidimensional reciprocity laws, \rm
 Math. Res. letters \bf 3 \rm (1996), 19-30

[{\bf Br-ML2}] J-L Brylinski and D. D. McLaughlin, \it 
 \v Cech cocycles for characteristic
classes,  Comm. Math. Phys. \bf 178 \rm (1996), 225-236

[{\bf Br-ML4}] J-L Brylinski and D. D. McLaughlin, \it 
 Non-commutative reciprocity laws
associated to finite groups, \rm Operads:
Proceedings of Renaissance Conferences, J-L. Loday,
J. D. Stasheff and A. A. Voronov eds, Contemp. Math. 
vol. \bf 202 \rm (1997), 421-438

[{\bf C-S}] J. Cheeger and J. Simons, {\it Differential characters
 and geometric invariants}, 
in Lecture Notes in Math. vol. 1167, Springer-Verlag (1985), pp. 50-80.

[{\bf Co}] A. Connes, \it Non-commutative differential geometry,
\rm Publ. Math. IHES \bf \rm (1986)

[{\bf E}] H. Esnault, \it Characteristic classes of flat bundles,
\rm Topology \bf 27 \rm (1988), 323-352

[{\bf E-V}] H. Esnault and E. Viehweg,
\it Deligne-Beilinson cohomology, 
\rm in
M.  Rapoport, N. Schappacher and 
P. Schneider, ed., 
Beilinson's Conjectures on Special 
Values of L-Functions,
\rm Perspectives in Math., \Ac ~(1988), 43-92

[{\bf Fe}] B. L. Feigin, \it On the cohomology of the Lie
algebra of vector fields and of the current algebra,
\rm Sel. Math. Sov. \bf 7 \rm (1988), 49-62

[{\bf Ka}] V. Kac, \it Infinite-Dimensional Lie algebras,
\rm Progress in Math. vol. \bf 44, \rm \Bi (1985)

[{\bf L-Q}] J-L. Loday and D. Quillen, \it
Cyclic homology and the Lie algebra homology
of matrices, \rm Comment. Math. Helv. \bf 59
\rm (1984), 565-591 

[{\bf Pi}] D. Pickrell, \it Extensions of loop groups, \rm
Algebras, Groups Geom. \bf 10 \rm (1993), 87-134

[{\bf P-S}] A. Pressley and G.Segal, {\it Loop Groups},
Clarendon Press (1986).

[{\bf Te}] C. Teleman, \it Some Hodge theory from Lie algebras, \rm
preprint (1998)

[{\bf Ts}] B. Tsygan, \it Homology of matrix Lie algebras
over rings and Hochschild homology, \rm Uspekh. Math. Nauk \bf 113
\rm (1983), 217-218

[{\bf W-X}] A. Weinstein and P. Xu, \it Extensions of symplectic
groupoids and quantization,
\rm J. Reine Angew. math. \bf 417 \rm (1991), 159-189
\vskip .22 in
Penn State University

Department of Mathematics

University Park, PA. 16802

e-mail: jlb@math.psu.edu

\bye